\documentclass{amsart}

\newtheorem{theorem}{Theorem}[section]
\newtheorem{proposition}[theorem]{Proposition}
\newtheorem{corollary}[theorem]{Corollary}
\newtheorem{lemma}[theorem]{Lemma}

\theoremstyle{definition}

\theoremstyle{remark}
\newtheorem{remark}[theorem]{Remark}

\numberwithin{equation}{section}

\def\CC{{\mathbb C}} \def\ZZ{{\mathbb Z}} 
 \def\QQ{{\mathbb Q}} 

\def\mal{\mathbin{\! \cdot \!}}
\def\Hom{\mathop{\rm Hom}\nolimits}
\def\t#1{\widetilde{#1}}
\def\b#1{\overline{#1}}
\def\Spec{\rm Spec}

\begin{document}

\title{On embeddings into toric prevarieties}

\author{J.~Hausen}
\address{Fachbereich Mathematik und Statistik, Universit\"at
Konstanz, D-78557 Konstanz, Germany.}
\email{Juergen.Hausen@uni-konstanz.de}

\author{S.~Schr\"oer}
\address{Mathematisches Institut, Ruhr-Universit\"at,
44780 Bochum, Germany.}
\email{s.schroeer@ruhr-uni-bochum.de}

\dedicatory{2 October 2002, final version}


\begin{abstract}
We give examples of  complete normal surfaces that are
not embeddable into
simplicial toric
prevarieties nor  toric prevarieties of affine intersection.
\end{abstract}

\maketitle

\section*{Introduction}

This note is concerned with embeddability and non-embeddability into
toric prevarieties. W\l odarczyk \cite{Wl}  has shown that a normal variety
$Y$ admits a closed embedding into a toric variety $X$ if and only if
every pair $y_1,y_2\in Y$ is contained in a common affine open
neighbourhood. If one drops the latter condition there still exists
a closed embedding into some toric {\it prevariety} $X$. This
means  that $X$ may be non-separated.

The above embedding result has the following refinement: Supposing
that $Y$ is $\QQ$-factorial one can choose $X$ to be simplicial and   
of affine intersection. 
In other words, $X$ is $\QQ$-factorial and the intersection of any two
affine open subsets of $X$ is again affine (see \cite{Wl} and
\cite{Ha}). 
Embeddings into toric prevarieties of affine intersection
can for example be used to obtain global resolutions of coherent
sheaves on $Y$ (see~\cite{Ha}).

The purpose of this note is to find out whether or not every
normal variety $Y$ can be embedded into a toric prevariety $X$ which
is simplicial or of affine intersection. Unfortunately, the
answer is negative: We provide  examples of normal  surfaces   that
admit  neither embeddings into toric prevarieties of affine
intersection nor into simplicial ones.

Together with W\l odarczyk's embedding result,   these
counterexamples imply that there are toric prevarieties of affine
intersection that cannot be embedded into simplicial toric
prevarieties of affine intersection. An analogous statement is easily seen
to hold  in the category of toric varieties, since there are
complete toric
varieties with trivial Picard group (see for example \cite{Ei}).

\section{Non-projective complete surfaces}

Here we briefly discuss some properties of complete normal
surfaces which are  non-projective. First, we have to fix   notation.
Throughout, we will work over an uncountable algebraically closed ground field
$k $.  The word \textit{prevariety} refers to an
integral
$k$-scheme of finite type. A \textit{variety} is a separated prevariety, and
a \textit{surface} is a 2-dimensional variety.
By a \textit{toric} prevariety we mean a prevariety
$X$ arising from a system of fans in some lattice $N$ in the sense of
\cite{acha}. For  $k=\CC$   this
notion precisely yields  the normal prevarieties $X$ endowed
with an effective regular action of a torus $T$ with an open
orbit.

The
following facts are well known (compare \cite{Ht}, theorem 4.2, and \cite{Wl}):

\begin{proposition}
Let $S$ be a normal surface and denote by $s_{1}, \ldots, s_{n}\in S$ its
non-$\QQ$-factorial singularities.
\begin{enumerate}
\item The surface $S$ is quasi-projective if and only if the points $s_{1},
   \ldots, s_{n}$ have a common affine neighbourhood.
\item The surface $S$ admits a closed embedding into a toric
   variety if and only if every pair $s_{i}, s_{j}$ lies in a common affine
   neighbourhood.
\end{enumerate}
\end{proposition}

In view of these statements, the simplest possible candidate for a
normal variety without closed embeddings into  toric prevarieties
of affine intersection is a normal surface having precisely two
non-$\QQ$-factorial singularities $s_{1}, s_{2}\in S$. In fact we will
work with surfaces of this type.

Let $S$ be a complete normal surface. There is a \textit{projective 
reduction}  $r \colon S
\to S^{\rm prj}$   of $S$, which means that  the morphism $r$ is universal
with respect to morphisms  to projective schemes (see \cite{BB}). Since
$\mathcal{O}_{S^{\rm prj}} \to r_*(\mathcal{O}_{S}) $ is necessarily 
bijective, the
projective scheme
$S^{\rm prj} $ is normal and the fibres of
$r $ are connected.

If $S^{\rm prj}$ is a curve,
any Weil divisor $E$ on $S$ has  a decomposition
$$ E = E^{\rm vert} + E^{\rm hor} $$
into the part $E^{\rm vert}$ consisting of those prime cycles of $E$
that are contained in the fibres of $r$ and the part $E^{\rm hor}$
which consists of the remaining prime cycles. A Weil divisor $E$ on
$S$ is called {\it vertical}  if $E = E^{\rm vert}$, and it is called
{\it horizontal} if $E = E^{\rm hor}$.

For a Weil divisor $D$ on a variety $Y$ and $y \in Y$, we denote by
$D_{y}$ the Weil  divisor obtained from $D$ by omitting all prime cycles
not containing the point $y$. The following statement will be
crucial for non-embeddability:

\begin{lemma}\label{vertical}
Let $S$ be a complete normal surface with precisely two non-$\QQ$-factorial
singularities $s_{1}$, $s_{2}\in S$. Assume that $S^{\rm prj}$ is
a curve and that every vertical Weil divisor on $S$ is
$\QQ$-Cartier. If $E$ is a $\QQ$-Cartier divisor on $S$ with $E_{s_{2}}
\ge 0$ and $E^{\rm hor} \ne 0$, then $E_{s_{1}}^{\rm hor}$ is not
effective.
\end{lemma}

\begin{proof}
Assume that $E_{s_{1}}^{\rm hor} \ge 0$ holds. Since $E$ and
$E^{\rm vert}$ are $\QQ$-Cartier,  so is their difference $E^{\rm hor}$.
Consider the decomposition
$$E^{\rm hor} = E^{\rm hor}_{+} + E^{\rm hor}_{-}$$
into positive and negative parts. Since $s_{1}$ and $s_{2}$ are not
contained in $E^{\rm hor}_{-}$, the surface $S$ is $\QQ$-factorial near
$E^{\rm hor}_{-}$. So  $E^{\rm hor}_{-}$ is $\QQ$-Cartier. Hence
$E^{\rm hor}_{+}$ is also $\QQ$-Cartier .
As $E^{\rm hor} \ne 0$, there must be at least one  horizontal Cartier divisor
$C\subset S$.

Consider the invertible sheaf $\mathcal{L}= \mathcal{O}_{S}(C)$.
Obviously,
$\mathcal{L}$ is $S^{\rm prj}$-generated on $S \setminus C$. Since $C
\to S^{\rm prj}$ is finite, we can apply \cite{Sc2},
Theorem 1.1 (over  affine open subsets
$\Spec(R)\subset S^{\rm prj}$), and deduce that $\mathcal{L}^{\otimes
n}$ is $S^{\rm prj}$-generated for some $n > 0$. Consequently, the
homogeneous spectrum
$$ S' = {\rm Proj}(f_{*} {\rm Sym}(\mathcal{L}))$$
is a projective $S^{\rm prj}$-scheme, hence projective. But
$\mathcal{L}$ is ample on the generic fibre of $r \colon S \to S^{\rm
prj}$,  contradicting the universal property of the projective
reduction.
\end{proof}

\begin{remark}\label{surfex}
Surfaces as above really exist, compare \cite{Sc}, 2.5. Here the assumption
that the ground field
$k $ is uncountable comes in.
\end{remark}

We will also make use of the following elementary fact:

\begin{lemma}\label{negativecomp}
Let $f \colon Y \to X$ be a morphism of integral normal
prevarieties. Given a Cartier divisor $D$ on $X$ such that the
preimage $E := f^{*}(D)$ exists as Cartier divisor. Decompose $D =
\sum_{i}
\lambda_{i} D_{i}$ and $E = \sum_{j} \mu_{j} E_{j}$  into prime cycles. If
there is a component $E_{j}$   with
$\mu_{j} < 0$, then there is a component $D_{i}$ with $\lambda_{i} <
0$ and $f(E_{j}) \subset D_{i}$.
\end{lemma}

\proof Assume there is no such $D_{i}$. Decompose $D=D_+-D_-$  into
positive and negative parts. Then we have
$E_{j} \not\subset f^{-1}(D_{-})$. Hence the restriction of $E$ to $Y
\setminus f^{-1}(D_{-})$ is not effective. On the other hand, the
restriction of $D$ to $X \setminus D_{-}$ is effective,
contradiction. \endproof

\section{Non-embeddability}

Recall that a scheme $X$ is separated if the diagonal morphism
$\Delta \colon X \to X \times X$ is a closed embedding. A weaker
condition is that $\Delta$ is affine. In this situation we call
$X$ to be of {\it affine intersection}. To check   this property it
suffices to find an open
affine covering $X = \bigcup_{i} X_{i}$  such that each
intersection $X_{i} \cap X_{j}$ is affine. Clearly, if $X$ is of 
affine intersection,
every subscheme  is so.

\begin{theorem}\label{sfmapsfactor}
   Let $S$ be a complete normal surface with precisely two
   non-$\QQ$-factorial singularities $s_{1},s_{2}\in S$. Assume
   that the projective reduction $S^{\rm prj}$ is a curve and that
   all irreducible  components of the fibres of $r \colon S \to S^{\rm prj}$ are
   $\QQ$-Cartier.  Let $f \colon S \to X$ be a morphism to a toric
   prevariety $X$. If $X$ is simplicial or of affine intersection, then
   there is a morphism  $\t{f} \colon   S^{\rm prj} \to X$ with $f =
   \t{f} \circ r$.
\end{theorem}

Let us record the following evident

\begin{corollary}\label{noembs}
   The surface $S$ is neither embeddable into  a toric prevariety of
   affine intersection nor into a simplicial toric prevariety.  \endproof
\end{corollary}

\noindent \textit{Proof of theorem \ref{sfmapsfactor}.}
First we treat the
case that $X$ is of affine intersection. Let $T$ denote the acting
torus of $X$. There is a unique $T$-orbit $B \subset X$
such that $f(S) \subset \b{B}$ holds. Note that $\b{B}$ is again a
toric prevariety of affine intersection. Replacing $X$ by $\b{B}$
we may assume that $f(S)$ hits the open $T$-orbit.

Set $x_{1} := f(s_{1})$ and $x_{2} := f(s_{2})$. For $i=1,2$ let
$X_{i} \subset X$ be the affine $T$-stable neighborhoods containing $T
\mal x_{i}$ as closed orbits, respectively. Set $N := \Hom(k^{*},T)$
and let $ M := \Hom(N,\ZZ)$. Then we have
$$
X_{i} = X_{\sigma_{i}}= {\rm Spec} (k[\sigma_i^\vee\cap M])
$$
for certain cones $\sigma_{i}$ in the lattice $N$. Since $X$ is of
affine intersection, the set $X_{12} := X_{1} \cap X_{2}$ is of the
form $X_{12} = X_{\sigma_{12}}$ for some common face $\sigma_{12}$ of
$\sigma_{1}$ and $\sigma_{2}$. Hence there is a linear form $u \in M
$ with
$$ u \in \sigma_{1}^{\vee}  \qquad {\rm and}\qquad \sigma_{12} =
\sigma_{1} \cap u^{\perp}. $$
Since $f(S)$ hits the
open $T$-orbit, all characters $\chi^{m}$, $m \in M$, admit pullbacks
$f^{*}(\chi^{m})$ as rational functions. In particular $\chi^{u}$ has
a pullback $\psi$ with respect to $f$, and the principal divisor ${\rm
   div}(\chi^{u}) $ has the principal divisor ${\rm div}(\psi)$ as
pullback. Since $\chi^{u}$ is defined on $X_{1}$, we have
$$ {\rm div}(\psi)_{s_{1}} \ge 0. $$

We claim that no component of
${\rm div}(\psi)_{s_{1}}$ contains the point $s_{2} \in S$. Otherwise,
let $E \subset S$ be such a component. Then there exists a $T$-stable
component $D_{1} \subset X_{1}$ of ${\rm div}(\chi^{u}) \cap {X_{1}}$
containing $f(E) \cap X_{1}$. This component corresponds to an edge
$\varrho_{1}\subset \sigma_{1}$. Note that by \cite{Fu}, p.~61 we
have
$\varrho_{1}
\not\subset u^{\perp}$. Take the closure $D := \b{D_{1}}$ in $X$,
set $D_{2} := D \cap X_{2}$, and let $\varrho_{2}\subset\sigma_{2}$ be
the edge corresponding to $D_{2}\subset X_2$. Then $s_{2} \in E$
implies
$$
x_{2} = f(s_{2}) \in f(E) \cap X_{2} \subset D_{2}. $$
Thus $D \cap
X_{2}$ is non-empty, hence $D_{12} := D \cap X_{12}$ is a
$T$-stable Weil divisor in $X_{12} = X_{\sigma_{12}}$. Let
$\varrho_{12}\subset\sigma_{12}$ be the edge  corresponding to
$D_{12}\subset X_{12}$. Since $D_{12}$ is induced by $D_{1}$ and
$D_{2}$ as well, we conclude $\varrho_{1} = \varrho_{12} =
\varrho_{2}$. Thus we have $\varrho_{1} \subset \sigma_{12} \subset
u^{\perp}$, contradicting $\varrho_{1} \not\subset u^{\perp}$.
As a consequence of our claim,  ${\rm div}(\psi)_{s_{1}}$ is an
effective $\QQ$-Cartier divisor. According to lemma~\ref{vertical},
the divisor ${\rm div}(\psi)_{s_{1}}$ is vertical.

Next, we show  that for every principal divisor $D = {\rm
   div}(\chi^{m})$, $m \in M$   its pullback to $S$ is vertical. It
   suffices to check this for  a set of generators of $M$, for example
   $M \cap  \sigma_{2}^{\vee}$. Assume that there is some $m \in M \cap
\sigma_{2}^{\vee}$ such that $E  := f^{*}(D)$ is not vertical. 
Decompose $D = \sum_{i}\lambda_i D_{i}$ and $E = \sum_{j}\mu_j E_{j}$ 
into prime
cycles. Since $D_{x_{2}} \ge 0$, we also have $E_{s_{2}} \ge 0$. By
lemma~\ref{vertical} there must be some horizontal component $E_{j}$
containing $s_1$ with $\mu_{j} < 0$. Lemma~\ref{negativecomp} tells us
that there is some component $D_{i}$ with $\lambda_{i} < 0$ and
$f(E_{j}) \subset D_{i}$. In particular we have $x_{1} \in D_{i}$.
Let $\varrho_{i}\subset \sigma_{1}$ be the edge corresponding to
the $T$-invariant Weil divisor $D_{i} \cap {X_{1}}$. Since $x_{2}\not\in
D_{i}$  we have $\varrho_{i} \not\subset \sigma_{12}$. Hence $\varrho_{i}
\not\subset \sigma_{1} \cap u^{\perp}$, and $D_{i}$ occurs with
positive multiplicity in ${\rm div}(\chi^{u})_{x_{1}}$.
On the other hand, we already have seen that ${\rm div}(\psi)_{s_{1}}$
is vertical. Consequently, $E_{j} \subset f^{-1}(D_{i})$ is also
vertical, contradiction. We have shown that ${\rm div}(f^{*}(
\chi^{m}))$ is vertical for all $m \in M$.

As a consequence, the induced morphism $f^{-1}(T) \to T$ is constant
along the generic fibre of $r \colon S \to S^{\rm prj}$. By the rigidity lemma
(see for example \cite{clkomo}, 1.5), all fibres of $r$ are mapped to 
points under $f
\colon S
\to X$. Now  \cite{EGA}, 8.11.1 gives a morphism $\t{f} \colon
S^{\rm prj}\to X$ with $f = \t{f} \circ r$. This proves the theorem for toric
prevarieties of affine intersection.

It remains to treat the case that $X$ is simplicial. According to
\cite{Ha}, Section~1, there is a toric prevariety $X'$ of affine
intersection and a local isomorphism $g \colon X \to X'$. As we have
seen, $g \circ f$ is constant along the fibres of $r \colon S \to
S^{\rm prj}$. Since the fibres of
$r $ are connected,  this implies that $f$ is constant along
the fibres of $r$, and we obtain the desired factorization as
above.
\endproof

\bibliographystyle{amsplain}

\end{document}